\documentclass[12pt]{article}
\usepackage{amsmath,amsthm}
\usepackage{amssymb,latexsym}
\usepackage{enumerate}
\usepackage[english]{babel}
\usepackage{graphicx}

\def\Q{{\mathbb Q}}
\def\Z{{\mathbb Z}}

\def\O{{\cal O}}

\newtheorem{lemma}{Lemma}
\newtheorem{theorem}[lemma]{Theorem}
\newtheorem{corollary}[lemma]{Corollary}

\newtheorem{proposition}[lemma]{Proposition}

\title{
Calculating power integral bases \\by using relative power integral bases\\
}
\author{
Istv\'{a}n Ga\'{a}l\thanks{
        Research supported in part by K115479 from the
        Hungarian National Foundation for Scientific Research
                         },\; 
L\'aszl\'o Remete and T\'\i mea Szab\'o
\\ \\
University of Debrecen, Mathematical Institute \\
H--4010 Debrecen Pf.12., Hungary \\
e--mail: igaal@science.unideb.hu, remetel42@gmail.com, \\ szabo.timea@science.unideb.hu
}

\begin{document}

\maketitle
\thispagestyle{empty}

\renewcommand{\thefootnote}{}



\renewcommand{\thefootnote}{\arabic{footnote}}
\setcounter{footnote}{0}

\begin{abstract}
Let $M\subset K$ be number fields.
We consider the relation of relative power integral bases of $K$ over
$M$ with absolute power integral bases of $K$
over $\Q$. We show how generators of absolute power integral bases 
can be calculated from generators of relative ones.
We apply our ideas in infinite families of octic fields
with quadratic subfields.
\end{abstract}

\section{Introduction: monogenity in the absolute and relative case}

Monogenity of number fields and the calculation of generators of power integral bases
is a classical topic of algebraic number theory c.f. \cite{nark}, \cite{book}. 
We have general algorithms for calculating generators of power integral bases in
lower degree number fields, \cite{gsch},
\cite{gppsim}, \cite{s5}, \cite{s6}. We only have partial results for higher degree
fields \cite{compos}, \cite{gop}, \cite{g6}, \cite{gp6i}, \cite{degnine}.

Let $K$ be an algebraic number field of degree $n$ with ring of integers $\Z_K$.
This field is {\it monogene} if $\Z_K$ is a simple ring extension of $\Z$, that is
there exist $\vartheta\in\Z_K$ such that $\Z_K=\Z[\vartheta]$. In this case
$\{1,\vartheta,\ldots,\vartheta^{n-1}\}$ is an integral basis of $K$, called
{\it power integral basis}. If $\alpha_1,\alpha_2\in\Z_K$ are related by  
$\alpha_1\pm \alpha_2\in\Z$ 
then the elements $\alpha_1,\alpha_2$ are called {\it equivalent}.
These elements have the same indices (see below) and $\alpha_1$ generates
a power integral basis of $K$ if and only if $\alpha_2$ does. Up to
equivalence there are only finitely many generators of power integral bases
of $K$.

We also considered monogenity and power integral bases  in the {\it relative case}
\cite{relcubic}, \cite{relquartic}, \cite{gsz2}.
The element $\vartheta$ generates a {\it relative power integral basis} 
of $K$ over the subfield $M$ if $\Z_K=\Z_M[\vartheta]$
($\Z_M$ denotes the ring of integers of $M$).
In the relative case we call $\alpha_1,\alpha_2\in\Z_K$ {\it equivalent} if 
$\alpha_1+\varepsilon \alpha_2\in \Z_M$ for some unit $\varepsilon$ in $M$.
These elements have the same relative indices (see below) and $\alpha_1$ generates
a relative power integral basis of $K$ over $M$ if and only if $\alpha_2$ does. Up to
equivalence there are only finitely many generators of relative power integral bases
of $K$ over $M$.

In the present paper we describe the relation of 
the generators of relative power integral bases
with the generators of absolute ones. 
We show how the generators of relative power integral bases
can be used to calculate generators of absolute power integral bases. 

The algorithm is especially simple if $M$ is a quadratic field.
We apply our method to three infinite families of octic fields with 
imaginary quadratic subfields.

\section{From relative power integral bases to absolute ones}
\label{relabs}

Let $M$ be an algebraic number field of degree $m$ and $K$ an extension of
$M$ with $[K:M]=k$. Then we have $[K:\Q]=k\cdot m$.
Let $\O$ be either the ring of integers $\Z_K$ of $K$ or an
order in $\Z_K$. Denote by $\Z_M$ the ring of integers of $M$.
We assume that there exist a relative integral basis of $\O$ over $M$.
(As we shall see in the following the existence of a power integral basis of $\O$
implies the existence of a relative power integral basis.)

Denote by $D_{\O}$ and $D_M$ the discriminants of $\O$ and the subfield $M$, respectively.
(In case $\O=\Z_K$ we have $D_{\O}=D_K$ where $D_K$ is the discriminant of
the field $K$.)
The {\it index} of a primitive element $\alpha$ of $\O$ with respect to the order $\O$ is
\begin{equation}
I_{\O}(\alpha)=\frac{\sqrt{|D(\alpha)|}}{\sqrt{|D_{\O}|}}.
\label{absind}
\end{equation}
We also have
\begin{equation}
I_{\O}(\alpha)=(\O^+ : \Z[\alpha]^+)=(\O^+ : \Z_M[\alpha]^+)\cdot 
(\Z_M[\alpha]^+: \Z[\alpha]^+),
\label{ind}
\end{equation}
where the indices of the additive groups of the corresponding rings are calculated.
The first factor is just the {\it relative index} of $\alpha$: 
\[
I_{\O/M}(\alpha)=(\O^+ : \Z_M[\alpha]^+).
\]
Denote by $D_{\O/M}$ the relative discriminant of ${\cal O}$ over $M$.
As it is well known
\begin{equation}
D_{\O}=N_{M/\Q}(D_{\O/M})\cdot D_M^{[K:M]}.
\label{discr}
\end{equation}
Denote by $\gamma^{(i)}$ the conjugates of any $\gamma\in M$
($i=1,\ldots,m$). Let $\delta^{(i,j)}$ 
be the images of $\delta\in K$ under the automorphisms of $K$ leaving the
conjugate field $M^{(i)}$ elementwise fixed ($j=1,\ldots,k$).
Then for any primitive element $\alpha\in \O$ we have
\[
I_{\O/M}(\alpha)=
\frac{\sqrt{|N_{M/\Q}(D_{\O/M}(\alpha))|}}{\sqrt{|N_{M/\Q}(D_{\O/M})|}}=
\]
\begin{equation}
=\frac{1}{\sqrt{|N_{M/\Q}(D_{\O/M})|}}
\cdot
\prod_{i=1}^m\;\;\prod_{1\leq j_1<j_2\leq k}\left|\alpha^{(i,j_1)}-\alpha^{(i,j_2)}\right|.
\label{relind}
\end{equation}
Further, by (\ref{absind}), (\ref{ind}), (\ref{discr})  and (\ref{relind})
we have
\[
J(\alpha)=(\Z_M[\alpha]^+: \Z[\alpha]^+)=
\]
\begin{equation}
=\frac{1}{\sqrt{|D_M|}^{[K:M]}}
\cdot
\prod_{1\leq i_1<i_2\leq m}\;\;\prod_{j_1=1}^k \prod_{j_2=1}^k
\left|\alpha^{(i_1,j_1)}-\alpha^{(i_2,j_2)}\right|.
\label{ind2}
\end{equation}

The element $\alpha$
generates a power integral basis of $\O$ if and only if $I_{\O}(\alpha)=1$.
Here we formulate the straightforward consequences of it, which will be very
useful in our calculations in the following sections.

By (\ref{ind}), $I_{\O}(\alpha)=1$ can only be satisfied if both 
factors of (\ref{ind}) are equal to 1.
Therefore, 
\begin{proposition}
A primitive element $\alpha\in\O$ generates a power integral basis of $\O$, 
if and only if 
\[ 
I_{\O/M}(\alpha)=1
\]
and
\begin{equation}
J(\alpha)=(\Z_M[\alpha]^+: \Z[\alpha]^+)=1.
\label{J}
\end{equation}
\end{proposition}
\noindent
Hence we have
\begin{corollary}
If $\alpha$ generates a power integral basis of $\O$, then it
generates a relative power integral basis of $\O$ over $M$.
\end{corollary}

It is well known that generators or relative power
integral bases are determined up equivalence, that is up 
to multiplication by a unit in $M$
and up to translation by element of $\Z_M$. Hence
\begin{proposition}
If $\alpha$ generates a power integral basis of $\O$, then
\begin{equation}
\alpha=A+\varepsilon \cdot \alpha_0,
\label{eps}
\end{equation}
where $\alpha_0$ is a generator of a relative power integral basis of $\O$ over $M$,
$\varepsilon$ is a unit in $M$ and $A\in\Z_M$.
\end{proposition}

Summarizing, in order to determine all generators of power integral bases of
$\O$ we have to perform the following steps:

\noindent
{\bf Step 1}
{\it Determine up to equivalence all generators $\alpha_0\in\O$ of 
relative power integral bases of $\O$ over $M$.}\\
In other words, determine all elements $\alpha_0\in\O$ with relative
index 1:
\[
I_{\O/M}(\alpha_0)=1.
\]
Note that if $\alpha_0$ has relative index 1, then by means of equivalence
any $\alpha$ of the form (\ref{eps}) also has relative index 1.

\noindent
{\bf Step 2}
{\it Given $\alpha_0$ determine $\varepsilon$ and $A$ so that $\alpha$ of (\ref{eps})
has $J(\alpha)=1$.}\\

\vspace{1cm}

Let $\mu_1=1,\mu_2,\ldots,\mu_m$ be an integral basis of $M$. Then the above $A$
can be represented in the form
\begin{equation}
A=a_1+a_2\mu_2+\ldots +a_k\mu_k.
\label{elem3}
\end{equation}
Since the (absolute) index is invariant under translation by an element
of $\Z$, we have to calculate $a_2,\ldots,a_m$ of (\ref{elem3}) up to sign.
Step 2 means to determine $\varepsilon$ and $a_2,\ldots,a_k$ satisfying
(\ref{J}). In view of (\ref{ind2}) this yields to solve an equation
of degree $k^2m(m-1)/2$ depending on $\varepsilon$ and $a_2,\ldots,a_k$.

This later task can became very complicated. However if $M$ is an imaginary
quadratic field, then there are only finitely many units $\varepsilon$ in $M$ and
we get a polynomial equation in one variable $a_2$. We shall apply our
method in this case in the following examples of infinite parametric
families of octic number fields.

\section{Simplest $D_4$ octics}

Recently B.K.Spearman and K.S.Williams \cite{d4} studied the family of simplest $D_4$
octics. Let $t>0$ be an integer parameter and $\vartheta$ a root of the polynomial
$x^8+(t^2+2)x^4+1$. They showed that these polynomials are irreducible, and the field
$K=\Q(\vartheta)$ has Galois group $D_4$. Assuming that $t^2+4$ is square free they 
calculated the discriminant of $K$ and gave an integral basis of $K$. By
\[
x^8+(t^2+2)x^4+1=(x^4+itx^2+1)(x^4-itx^2+1)
\]
$M=\Q(i)$ is a subfield of $K$.

Here we restrict ourselves to parameters of the form $t=2T^2$.
We explicitely describe all generators of relative power integral bases of 
the order ${\cal O}=\Z_M[\vartheta]$ over $\Z_M$. 
Moreover we show that the order $\cal O$ admits no power integral bases.

\subsection{Relative power integral bases in the family of $D_4$ octics}

Let $T$ be a nonzero integer parameter and
$K$ the algebraic number field generated by a root $\vartheta$
of the polynomial $f(x)=x^8+(4T^4+2)x^4+1$. 
Let $M=\Q(i)$.  Denote by $\Z_K$ (resp. $\Z_M$) the ring of integers of
$K$ (resp. $M$). 
Consider the order ${\cal O}=\Z_M[\vartheta]$ of $K$.

Our purpose is to explicitely determine all generators of relative
power integral bases of $\cal O$ over $M$. Obviously, any $\alpha\in\cal O$
can be written in the form
\begin{equation}
\alpha=A+X\vartheta+Y\vartheta^2+Z\vartheta^3
\label{alfa}
\end{equation}
with $A,X,Y,Z\in\Z_M$.

\begin{theorem}
Assume $T>11$.
Up to equivalence all generators of relative power integral bases of $\cal O$ over $\Z_M$
are given by
\[
\alpha=\vartheta, 
\]
\[
\alpha=-2iT^2\vartheta+\vartheta^3, 
\]
\[
\alpha=(1+4T^4)\vartheta\pm (1+i)T\vartheta^2+2iT^2\vartheta^3,
\]
\[
\alpha=\pm(1+i)T\vartheta^2+\vartheta^3.
\]
\label{th1}
\end{theorem}

\noindent
{\bf Proof of Theorem \ref{th1}}.
The octic polynomial $f(x)$ can be written as
\[
f(x)=(x^4+2iT^2x^2+1)(x^4-2iT^2x^2+1)
\]
hence the relative defining polynomial of $\vartheta$ over $M$
is $x^4-2iT^2x^2+1$.
In our proof we use the result of I.Ga\'al and M.Pohst \cite{relquartic} 
on power integral bases in relative quartic extensions (cf. also \cite{book}).

According to \cite{relquartic} the coefficients $X,Y,Z\in\Z_M$ of $\alpha$
in (\ref{alfa}) must satisfy
\begin{eqnarray*}
F(U,V)=(U-2iT^2V)(U-2V)(U+2V)&=&\varepsilon,\\
Q_1(X,Y,Z)=X^2-2iT^2Y^2+4iT^2XZ+(1-4T^4)Z^2&=&U,\\
Q_2(X,Y,Z)=Y^2-XZ-2iT^2Z^2&=&V,
\end{eqnarray*}
with a unit $\varepsilon$ of $M$ and with $U,V\in\Z_M$. 
We have to determine the solutions $U,V\in\Z_M$ of the first equation and for all 
pairs $U,V$ to calculate the corresponding solutions $X,Y,Z$ of the second and third
equations. By the first equation we have
$U-2V=\varepsilon_1$ and $U+2V=\varepsilon_2$ with units $\varepsilon_1,\varepsilon_2\in M$.
Therefore $4V=\varepsilon_2-\varepsilon_1$. Since all units in $M$ are $\pm 1, \pm i$,
the only $V\in\Z_M$ satisfying this equation is $V=0$. Hence $U$ is again
a unit in $M$. Following the method of \cite{relquartic} we set
\[
Q_0(X,Y,Z)=UQ_2(X,Y,Z)-VQ_1(X,Y,Z)=0.
\]
Using standard arguments described in \cite{relquartic} we can parametrize
$X,Y,Z$ with parameters $P,Q\in\Z_M$ so that up to a unit factor we get
\begin{equation}
X=P^2-2iT^2Q^2,\;\; Y=PQ,\;\; Z=Q^2.
\label{kxyz}
\end{equation}
Substituting the formulas (\ref{kxyz}) into $Q_1(X,Y,Z)=U$ we obtain a quartic
relative Thue equation over $M$:
\begin{equation}
P^4-2iT^2P^2Q^2+Q^4=\varepsilon,
\label{relthue1}
\end{equation}
with a unit $\varepsilon$ in $M$. This equation can be written in the
form
\begin{equation}
P^4-((1+i)T)^2P^2Q^2+Q^4=\varepsilon,
\label{relthue2}
\end{equation}
therefore we may apply the results of V.Zielger \cite{ziegler} on the solution 
of this equation by taking $t=(1+i)T$ as parameter.
Theorem 2 of \cite{ziegler} implies that,
assuming $|t^2|>245$, that is $|T|>11$, 
up to unit factors of $M$ all solutions of (\ref{relthue1}) are
\begin{equation}
(P,Q)=(1,0),\; (0,1), \; (1,\pm (1+i)T), \; ((1+i)T,\pm 1).
\label{pq}
\end{equation}
Substituting these vales of $(P,Q)$ into (\ref{kxyz}) we obtain the
possible triplets:
\begin{equation}
\begin{array}{|c|c|c|c|}
\hline
&x&y&z\\\hline
Case \; 1 & 1  &  0&  0 \\ \hline
Case \; 2 & -2iT^2 &  0&  1 \\ \hline
Case \; 3 & 1+4T^4  &  \pm (1+i)T &  2iT^2 \\ \hline
Case \; 4 & 0 &  \pm (1+i)T &  1 \\ \hline
\end{array}
\label{xxxx}
\end{equation}
This proves Theorem \ref{th1}. $\Box$

\subsection{Power integral bases in the family of $D_4$ octics}

Despite of the promising result on relative power integral bases we have

\begin{theorem}
For $|T|>11$ the order $\cal O$ admits no power integral bases.
\label{th2}
\end{theorem}

\noindent
{\bf Proof of Theorem \ref{th2}}.
In view of (\ref{eps}) a generator $\alpha$ of a power integral basis of $\O$
must be of the form $\alpha=a_1+a_2i+\varepsilon \alpha_0$ where 
$a_1,a_2\in\Z$, $\varepsilon=\pm 1,\pm i$ and the
possibe values of $\alpha_0$ are listed in Theorem \ref{th1}.
Any $\alpha$ of the above form has relative index $I_{\O/M}(\alpha)=1$.
The index of $\alpha$ is independent of $a_1$ and it is sufficient to determine
$\alpha$ up to sign. Therefore we have to consider the possible 
values of $\alpha_0$ and for $\varepsilon=1$, $\varepsilon=i$ and we have to calculate 
$J(\alpha)$. We have $D_M=-4$ hence
\[
J(\alpha)=\frac{1}{2^4}\prod_{j_1=1}^4 \prod_{j_2=1}^4 
\left| \alpha^{(1,j_1)}- \alpha^{(2,j_2)}\right|.
\]
In Case 1 we get
\[
J(\alpha)=2^4\cdot |(4T^2a_2^2-1+4a_2^2)(4T^2a_2^2+1-4a_2^2)
(T^8+8a_2^4T^4+16a_2^8+16a_2^4)|.
\]
Hence $J(\alpha)$ is divisible by $2^4$, yielding that
$\alpha$ can not be a generator of a power integral basis.

In the other cases we got much more complicated formulas, but in each case
$J(\alpha)$ is divisible by $2^4$. $\Box$

\subsection{Remarks on the numerical calculations}

All calculations involved in the proof of Theorem \ref{th2} 
were performed in Maple \cite{maple} under Linux. 
$J(\alpha)$ is a polynomial with integer coefficients
of degree 16 in $a_2$, depending also heavily on the parameter $T$. 
We used symmetric polynomials several times to simplify the formulas.
Without being very careful the formulas became extremely complicated
and Maple broke down in lack of memory space. Using careful approach
all calculations took less than 2 minutes.

\section{Composites of imaginary quadratic fields and
pure quartic fields}

In a recent paper \cite{grsz} we considered number fields of 
type $K=\Q(\sqrt[4]{m},i\sqrt{d})$ for $d=3,7,11,19,43,67,163$
and for $1<m\leq 5000$, $m \equiv 2,3 (\bmod \ 4)$ with $(d,m)=1$.
Set $\xi=\sqrt[4]{m},\omega=(1+i\sqrt{d})/2$, then 
\[
\{1,\xi, \xi^2, \xi^3,\omega,\omega\xi, \omega\xi^2, \omega\xi^3\}
\]
is an integral basis of $K$ and $\{1,\xi, \xi^2, \xi^3\}$ is a relative 
integral basis of $K$ over $M=\Q(i\sqrt{d})$.  In \cite{grsz}
we described all generators
\[
\alpha=A+X\xi+Y\xi^2+Z\xi^3
\]
of relative power integral bases of $K$ over $M$ with 
$A,X,Y,Z\in\Z_M$ and
$\max(\overline{|X|},\overline{|Y|},\overline{|Z|})<10^{500}$
(here $\overline{|X|}$ denotes the size of $X$ that is the
maximum absolute value of its conjugates).
The problem lead us to a quartic relative binomial Thue equation.
Using the algorithm of \cite{girelthue} we calculated the "small" solutions
of this equation which resulted Theorem 3 of \cite{grsz}. Note that 
according to our experience these equations never have "large" solutions
hence our list contains all solutions with high probability.
Further, calculating the "small" solutions was the only way to deal
with thousands of relative Thue equations.

Using the ideas of Section \ref{relabs} we tested if 
there exist generators of power integral bases of $K$ over $\Q$
corresponding to the relative power integral bases found in Theorem 3 
of \cite{grsz}.
We have
\begin{theorem}
Let $d=3,7,11,19,43,67,163$ and $1<m\leq 5000$ with $m \equiv 2,3 (\bmod \ 4)$ 
and $(d,m)=1$. Then the number field $K=\Q(\sqrt[4]{m},i\sqrt{d})$
does not admit any generators of power integral bases of the form
\[
\alpha=A+\varepsilon(X\xi+Y\xi^2+Z\xi^3)
\]
where $A\in \Z_M$, $\varepsilon$ a unit in $M$ and 
$X,Y,Z\in \Z_M$ with
\[
\max(\overline{|X|},\overline{|Y|},\overline{|Z|})<10^{500}.
\]
\label{dm}
\end{theorem}

\noindent
{\bf Proof of Theorem \ref{dm}}.
For all possible values of $X,Y,Z$ listed in Theorem 3 of \cite{grsz}
and for all possible unit $\varepsilon$ in $M$ we 
set $A=a_1+a_2\omega$. We calculated $J(\alpha)$ which is a polynomial 
in $a_2$ with integral coefficients of degree 16. In each case we found
that $J(\alpha)=\pm 1$ is not solvable for $a_2$ in integers.
Calculation with polynomials with integer coefficients was very fast,
the whole calculation took a few seconds.  $\Box$

\section{Parametric families of quartic extensions of imaginary quadratic fields}

In \cite{gsz2} we calculated generators of relative power integral bases in 
infinite parametric families of orders of certain octic fields.
Here in two of these families we check  if there exist corresponding
generators of (absolute) power integral bases.
The challenge of these examples is that $J(\alpha)$ depends not only on $a_2$
but also on the quadratic field and the parameter of the family.

\vspace{1cm}

{\bf I.} Let $d>0$ be an integer, $-d\equiv 2,3 (\bmod{4})$ and set $M=\Q(i\sqrt{d})$.
Let $t\in\Z_M$ be a parameter and let $\xi$ be a root of the polynomial
\[
f(x)=x^4-t^2x^2+1.
\]
Consider $\O=\Z_M[\xi]$. In \cite{gsz2} we showed that 
for $|t|>245$ up to equivalence there are 
five generators of relative power integral bases of $\O$ over $\Z_M$, namely
\[
\alpha_0=\xi,-t^2\xi+\xi^3,(1-t^4)\xi+t\xi^2+t^2\xi^3,(1-t^4)\xi-t\xi^2+t^2\xi^3,
t\xi^2+\xi^3,-t\xi^2+\xi^3.
\]
We have
\begin{theorem}
Under the above conditions for $|t|>245$ the order $\O$ admits no power integral bases.
\label{dm1}
\end{theorem}

\noindent
{\bf Proof of Theorem \ref{dm1}}.
Denote by $\alpha_0$ a possible generator of a relative power integral basis
of $\O$ over $\Z_M$, say
\[
\alpha_0=(1-t^4)\xi+t\xi^2+t^2\xi^3
\]
where $t=t_1+t_2i\sqrt{d}$ is the parameter ($t_1,t_2\in\Z$).
Note that since the minimal polynomial of $\xi$ over $\Z_M$
depends on the parameter $t\in\Z_M$, hence 
$\xi$ depends on $t$ and also on $d$. We let $\varepsilon=\pm 1$ and represent
$\alpha$ in the form
\[
\alpha=a_1+a_2i\sqrt{d}+\varepsilon \alpha_0
\]
with $a_1,a_2\in\Z$.
Then we calculate $J(\alpha)$. This is a very complicated polynomial of degree 16
depending not only on $a_2$ but also on $t_1,t_2,d$. Using symmetric polynomials
and simplifying the formulas very carefully, we obtain that $J(\alpha)$ is 
disvisible by 16. Therefore there are no generators of power integral bases of 
$\O$ corresponding to $\alpha_0$. The proof runs the same way for the other four candidates
of $\alpha_0$, as well. The Maple calculation took 10-60 seconds per case. $\Box$

\vspace{1cm}

{\bf II.} Let again $d>0$ be an integer, $-d\equiv 2,3 (\bmod{4})$, $M=\Q(i\sqrt{d})$.
Let $t\in\Z_M$ be a parameter and let $\xi$ be a root of the polynomial
\[
f(x)=x^4-4tx^3+(6t+2)x^2+4tx+1.
\]
Let $\O=\Z_M[\xi]$. According to \cite{gsz2} for
$|t|>1544803$ up to equivalence there are 
two generators of power integral bases of $\O$ over $\Z_M$, namely
\[
\alpha_0=\xi,(6t+2)\xi-4t\xi^2+\xi^3.
\]
We have
\begin{theorem}
Under the above conditions for $|t|>1544803$ the order $\O$ admits no power integral bases.
\label{dm2}
\end{theorem}

The proof of this statement is similar to the proof of Theorem \ref{dm1}.

\end{document}